\documentclass[11pt]{amsart}
\usepackage[marginratio=1:1,height=615pt,width=360pt,tmargin=117pt]{geometry}
\usepackage{amssymb,amsthm,amsmath}
\usepackage[numbers,sort&compress]{natbib}
\usepackage{color}
\usepackage{graphicx}
\usepackage{float}
\usepackage{tikz}
\usepackage[colorlinks,
            linkcolor=blue,
            anchorcolor=blue,
            citecolor=blue
            ]{hyperref}

\def\e{\varepsilon}

\let\newpf\proof \let\proof\relax 
\newenvironment{pf}{\newpf[\proofname]}{\qed\endtrivlist}

\newcommand{\ba}{\overline{A}}

\def\be{\begin{equation}}
\def\ee{\end{equation}}

\def\ba{{\begin{align}}}
\def\ea{{\end{align}}}

\def\bm{\begin{matrix}}
\def\em{\end{matrix}}

\def\0{{\mathbf 0}}

\newtheorem{Theorem}{Theorem}[section]
\newtheorem{Lemma}{Lemma}[section]
\newtheorem{Proposition}{Proposition}[section]

\theoremstyle{remark}
\newtheorem{Remark}{Remark}[section]

\newtheorem{Definition}{Definition}[section]

\numberwithin{equation}{section}

\theoremstyle{definition}

\newcommand{\dist}{\operatorname{dist}}

\newcommand{\Q}{{\mathbb Q}}
\newcommand{\R}{{\mathbb R}}
\newcommand{\T}{{\mathbb T}}

\newcommand{\Z}{{\mathbb Z}}

\def\B0{{\bold{0}}}

\catcode`\@=12

\def\Empty{}
\newcommand\oplabel[1]{
  \def\OpArg{#1} \ifx \OpArg\Empty {} \else
    \label{#1}
  \fi}

\newcommand{\comm}[1]{}
\newcommand{\comment}[1]{}

\begin{document}

\title[]{Arithmetic version of anderson localization for quasiperiodic Schr\"odinger operators with even cosine type potentials}
\author{}
\author{Lingrui Ge}
\address{
Department of Mathematics, University of California Irvine, CA, 92697-3875, USA}

 \email{lingruige10@163.com}

\author {Jiangong You}
\address{
Chern Institute of Mathematics and LPMC, Nankai University, Tianjin 300071, China} \email{jyou@nankai.edu.cn}

\author{Xin Zhao}
\address{
Department of Mathematics, University of California Irvine, CA, 92697-3875, USA}

 \email{njuzhaox@126.com}

\begin{abstract}
We propose a new method to prove Anderson localization for quasiperiodic Schr\"odinger operators and apply it to the quasiperiodic model considered by Sinai \cite{Sinai} and  Fr\"ohlich-Spencer-Wittwer \cite{fsw}. More concretely, we prove Anderson localization for even $C^2$ cosine type quasiperiodic Schr\"odinger operators with large coupling constants, Diophantine frequencies and Diophantine phases.
\end{abstract}
\maketitle
\section{Introduction}

Anderson localization of particles and waves in disordered media is one of the most intriguing phenomena in solid-state physics found by Anderson \cite{anderson}. Mathematically, localization  has already been extensively studied for the random  cases \cite{a,am,asfh,dks,ks}.  Anderson localization in quasiperiodic media even has stronger  backgrounds in physics \cite{aos,oa}, which is known to have strong connection with integer quantum Hall effect \cite{ag,h,ntw}, and also plays an important  role  in the emerging subject of optical crystals \cite{fk}.

Although been studied for over sixty years, Anderson localization \footnote{Pure point spectrum with exponentially decaying eigenfunctions} for quasiperiodic operators has not been completely understood since it depends sensitively on the frequency, the phase and amplitude of oscillation of the potential.
So far almost sure Anderson localization \footnote{Anderson localization for almost every phase.} for {\it fixed} Diophantine frequencies was rigorously proved only for  the following {\it cosine type} quasiperiodic Schr\"odinger operators with large coupling constants by Sinai \cite{Sinai}, Fr\"ohlich-Spencer-Wittwer \cite{fsw} \footnote{Fr\"ohlich-Spencer-Wittwer \cite{fsw} requires the potential to be even.} and Forman-Vandenboom \cite{fv}
\begin{equation}\label{1.1}
	(H_{\lambda v,\alpha,\theta}u)_n=u_{n+1}+u_{n-1}+\lambda v(\theta+n\alpha)u_n, \forall n\in\Z,
\end{equation}
with $\alpha\in\R\backslash\Q$ (the frequency), $\theta\in\T:=\R/\Z$ (the phase), $\lambda\in\R$ (the coupling constant) and $v\in C^2(\T,\R)$ (the potential) satisfying
\begin{itemize}
	\item $\frac{dv}{d\theta}=0$ at exactly two points, one is minimal and the other one is maximal, which are denoted by $z_1$ and $z_2$.
	\item These two extremals are non-degenerate, that is, $\frac{d^2v}{d\theta^2}(z_j)\neq0$ for $j=1,2.$
\end{itemize}

 Anderson localization (AL)  with precise arithmetic descriptions on  both the frequencies and full Lebesgue measure set of the phases, is referred as to  {\it arithmetic version of Anderson localization}, which is obviously stronger than almost sure Anderson localization.  The first arithmetic version of Anderson localization was given by Jitomirskaya \cite{J} who proved  that for any {\it fixed} Diophantine frequency and {\it any  fixed} Diophantine phase, the almost Mathieu operator \footnote{Operator \eqref{1.1} with $v(\theta)=2\cos2\pi\theta$.} has AL for $|\lambda|>1$.  Such arithmetic description on the frequency and the phase was explored in a sharp way by Jitomirskaya and Liu, namely, for Diophantine phase, there is a sharp spectral transition in frequency \cite{JLiu}, and for Diophantine frequency, there is a sharp spectral transition in phase \cite{JLiu1}. Arithmetic Anderson localization for one dimensional long range quasiperiodic operators with cosine potential was proved in  \cite{bj02,aj1}. Recently, a new nonperturbative proof of arithmetic theoretic Anderson localization was given in \cite{gy}, which applies to higher dimensional long range quasiperiodic operators,  based on nonperturbative reducibility method and duality argument.

  However, the proofs of all the above arithmetic Anderson localization results crucially depend on the assumption that the potential is exactly the {\it cosine} function. It is not clear if arithmetic Anderson localization could be expected for other potentials.
The main purpose of this paper is to present a new method from the point of view of dynamical systems  to prove the arithmetic version of Anderson localization for quasiperiodic Schr\"odinger operators. Applying our method to even cosine type quasiperiodic Schr\"odinger operator, we give an improvement of Fr\"ohlich-Spencer-Wittwer's result \cite{fsw}.  Compared to the methods of \cite{Sinai,fsw,fv} which are based on certain kind of multiscale analysis, our method is purely dynamical and gives concrete description of the localization phases \footnote{In \cite{Sinai,fsw,fv}, Anderson localization was proved for almost every phase without an arithmetic description.}. We are also able to give almost sharp estimate on the decay rate of the all eigenfunctions.
\subsection{Statement of the main results}
Before formulating our results, we first give precise arithmetic description on $\alpha$ and $\theta$. A frequency  $\alpha\in \R$ is called ($\kappa,\tau$)-{\it Diophantine}   (denoted by $\alpha\in \rm{DC}(\kappa,\tau))$ if
\begin{equation}\label{dc_def}
\dist(k\alpha,\Z)\ge \gamma|k|^{-\tau},\quad \forall k\in \Z\backslash\{0\}.
\end{equation}
We will use the notation
$$
\rm{DC}=\bigcup_{\kappa>0;\,\tau>1} \rm{DC}(\kappa,\tau).
$$
For a given irrational number $\alpha$, we say $\theta\in (0,1)$ is ($\gamma,\tau$)--Diophantine with respect to  $\alpha$ (denoted by $\Theta_\gamma^\tau$) if
$$
\|2\theta+k\alpha\|_{\R/\Z}>\frac{\gamma}{(|k|+1)^\tau},
$$
for any $k\in\Z$, where $\|x\|_{\R/\Z}=\text{dist}(x,\Z).$ Let $
\Theta=\bigcup_{\gamma>0;\,\tau>1}\Theta_\gamma^\tau.
$
 Clearly, $\Theta$ is a set of  full Lebesgue measure for any fixed irrational number $\alpha$.

\begin{Theorem}\label{main}
Given $\alpha\in DC$ and an even $C^2$ cosine type potential $v$, there exists  $\lambda_0(\alpha,v)$ such that $H_{\lambda v,\alpha,\theta}$ has Anderson localization for all $\theta\in \Theta$ provided that $\lambda>\lambda_0$.
\end{Theorem}
\begin{Remark}
To give a simple arithmetic description of the localization phases (i.e., the Diophantine phases), the eveness condition seems to be necessary.
\end{Remark}
\begin{Remark}
If $\alpha$ is very Liouvillean or $\theta$ is very $\alpha$-Liouvillean (i.e., for generic $\alpha$ and $\theta$), $H_{\lambda v,\alpha,\theta}$ has purely singular continuous spectrum \cite{ayz,JLiu1,js}. Thus to prove localization type results, the arithmetic assumptions on both $\alpha$ and $\theta$ are necessary.
\end{Remark}
We also have a precise estimate on the decay rate of all eigenfunctions.
\begin{Theorem}
Given $\alpha\in DC$, $\e>0$ and an even $C^2$ cosine type potential $v$, there exists $\lambda_0(\alpha,v,\e)$ such that  all eigenfunctions of the operator $H_{\lambda v,\alpha,\theta}$  satisfy
$$
\liminf\limits_{|n|\rightarrow\infty}-\frac{\ln\left(u^2_E(n)+u^2_E(n+1)\right)}{2|n|}\geq (1-\e)\ln\lambda
$$ provided that $\lambda>\lambda_0$.
\end{Theorem}
\begin{Remark}
For the almost Mathieu operator (a typical example), Jitomirskaya \cite{J} proved
$$
\liminf\limits_{|n|\rightarrow\infty}-\frac{\ln\left(u^2_E(n)+u^2_E(n+1)\right)}{2|n|}=\ln\lambda.
$$
Thus, the decay rate in the above theorem is almost sharp.
\end{Remark}
We point out an interesting phenomenon based on Theorem \ref{main}: The localization phases do not  sensitively depend on the space of even $C^2$ cosine type $v$. This phenomenon can be viewed as the robustness of localization phases introduced in \cite{gy}.
\begin{Definition}
For fixed $\alpha$, $H_{V,\alpha,\theta}$ is said to have $C^r$ robust Anderson localization if there is a $C^r$ neighborhood $B(V)$ of $V$ and a subset $\tilde{\Theta}$, such that
$$
\bigcap_{\tilde V \in B(V)}\{ \theta \, | H_{\tilde{V},\alpha,\theta}\ \  has \ \ AL\}=\tilde\Theta,
$$
moreover $|\tilde \Theta|=1$.
\end{Definition}
Theorem \ref{main} proved that $H_{\lambda v,\alpha,\theta}$ with even cosine like potential has $C^2$ robust Anderson localization in the space of even potentials. It seems that  both the symmetry (eveness) and the profile ($C^2$ cosine type) of the potential play key roles in robust Anderson localization.

We also mention some important results related to Anderson localization.  Eliasson \cite{Eli97}  proved that if $v$ is a Gevrey function satisfying non-degenerate conditions, for any fixed Diophantine $\alpha$, $H_{\lambda v,\alpha,\theta}$ has pure point spectrum for $a.e.$ $\theta$ and large enough $\lambda$ (depending on $\alpha$).  Bourgain and Goldstein \cite{bg} proved that, in the positive Lyapunov exponent regime, for any  fixed $\theta$, $H_{\lambda v,\alpha,\theta}$ has AL for  $a.e.$  Diophantine $\alpha$ provided that  $v$ is a non-constant real analytic function. See \cite{bj,K,gyz,gk} for more results.

\subsection{Strategy of the proof} As we mentioned above, our method is motivated by the methods introduced in \cite{ayz,gy,wz1,gyzh1,jk2}. In \cite{ayz,jk2}, Avila-You-Zhou and Jitomirskaya-Kachkovskiy gave criteria to prove almost sure Anderson localization for quasiperiodic operators based on nonperturbative reducibility method and duality argument. More precisely, there are two steps. The first step is to construct a family of eigenvalues and eigenfunctions for almost every phase  based on reducibility and duality. The second step is to show the family of eigenfunctions they constructed form a complete basis of $\ell^2(\Z)$. However, the  arithmetic version of  Anderson localization is more difficult to be proved compared with almost sure Anderson localization. In \cite{gy}, Ge-You found a strategy to recover the phases lost in using the method in \cite{ayz,jk2}, by introducing an auxiliary measure defined by reducibility, i.e. the $\mathcal{R}$-measure. By using quantitative reducibility, they proved stratified continuity of the $\mathcal{R}$-measure with respect to the phase on the set of Diophantine phases. Ge-You's method was further developed and simplified by Ge-You-Zhao in \cite{gyzh1} to give a new proof of the arithmetic transition conjecture proposed by Jitomirskaya \cite{j}.

In the present paper, instead of using reducibility and duality, we give a new way to construct a family of eigenvalues with exponentially localized eigenfunctions for $C^2$ cosine type quasiperiodic Schr\"odinger operators by the induction scheme developed in \cite{wz1}. The intuition is that if the intersection between asymptotic stable and unstable directions of the transfer matrix persist in larger and larger time scale which eventually implies the intersection of stable and unstable directions and  the norm of the transfer matrix grows exponentially, then one can construct an eigenfunction. When a family of eigenvalues with exponentially localized  eigenfunctions are constructed almost surely, almost sure Anderson localization follows directly by the criteria in \cite{ayz,jk2}.  To prove  the arithmetic version of AL, i.e., AL for all $\theta\in \Theta$,   one possible way is to prove that  $d\mu^{pp}_\theta$ \footnote{the pure point piece of the spectral measure} is continuous in $\Theta$. However, this seems to be a difficult task and we don't know how to prove it directly since $d\mu^{pp}_\theta$ sensitively depends on $\theta$.   Our strategy is to  introduce a new measure $d\nu_\theta$ via the localized eigenfunctions we constructed, which is called  $\mathcal{L}$-measure, motivated by the $\mathcal{R}$-measure defined in \cite{gy}. We will prove that  $d\nu_\theta$ is absolutely continuous with respect to $d\mu^{pp}_\theta$. The advantage of $d\nu_\theta$ is its stratified continuity in $\Theta$, more precisely the continuity in $\Theta_\gamma^\tau$, can be proved by  quantitative estimates of the localized eigenfunctions. In this way, we can approximate each lost phase in  $ \Theta$ by localization phases, and prove  $d\mu_\theta^{pp}(\R)= d\nu_\theta(\R)=1$ for all phases in $ \Theta$.

\section{ Wang-Zhang's induction theorem}
Inspired by \cite{young}, Wang and Zhang in \cite{wz1} developed an induction scheme to study the positivity and continuity of the Lyapunov exponent and Cantor spectrum of quasiperiodic Schr\"odinger operators  with   $C^2$ cosine type potentials \cite{wz1,wz2}.   We observe that the induction theorem can also be used to construct eigenfunctions, which is one of the corner stones in our proof of the arithmetic version of Anderson localization. Now we briefly introduce their induction theorem.  The readers are referred to \cite{wz1} for details.

For $\theta \in \R/\Z,$ let
$$
R_{\theta}=\begin{pmatrix}
\cos{2\pi\theta}&-\sin{2\pi\theta}\\
\sin{2\pi\theta}&\cos{2\pi\theta}
\end{pmatrix}
\in SO(2,\R).
$$
Define
$$s:SL(2,\R)\rightarrow\R{\mathbb{P}}^1=\R/(\pi\Z)$$
as  the most contraction direction of $A\in SL(2,\R),$ i.e.,  $\|A\cdot\hat{s}(A)\|=\|A\|^{-1} $ for unit vector $\hat{s}(A)$ in the direction $  s(A)$. Abusing the notation a little, define $u(A)=s(A^{-1})$ and $\hat{u}(A)= \hat s(A^{-1})$. Then for $A\in SL(2,\R)$, it is clear that
$$
A=R_u\cdot\begin{pmatrix}
\|A\|&0\\
0&\|A\|^{-1}
\end{pmatrix}
\cdot R_{\frac{\pi}{2}-s},
$$
where $s,u\in[0,\pi)$ are angles corresponding to the directions $s(A),u(A)\in \R/(\pi\Z).$

The main object of the present paper is the following \textit{Schr\"{o}dinger cocycles} $(\alpha,S_E^{\lambda v})$, where
$$
S_E^{\lambda v}(x):=
\begin{pmatrix}
E-\lambda v(x) & -1\\
1 & 0
\end{pmatrix},   \quad E\in\R.
$$
Let $\lambda\ge \lambda_0=\lambda_0(v)\gg 1$ and $t=\frac{E}{\lambda}\in J=[\inf v-2, \sup v+2]$. In this case,  there is $B\in C^2(\T\times J,SL(2,\R))$ such that
\begin{equation}\label{schcon}
B^{-1}(x+\alpha,t)\begin{pmatrix}E-\lambda v(x)&-1\\ 1&0
\end{pmatrix}B(x,t)=A(x,t),
\end{equation}
where
\begin{align*}
\nonumber A(x,t)=\begin{pmatrix}
\lambda(x,t)& 0\\
0& \lambda^{-1}(x,t)
\end{pmatrix}\cdot R_{\phi(x,t)},\ \ \cot\phi(x,t)=t-v(x).
\end{align*}
\eqref{schcon} in fact  gives the polar decomposition of the Schr\"odinger cocycles.\\

From now on, let $A(x,t)$ be as above and
$$
A_n(x,t):=
\left\{\begin{array}{l l}
A(x+(n-1)\alpha,t) \cdots A(x+\alpha,t) A(x,t),  & n\geq 0\\[1mm]
A^{-1}(x+n\alpha,t) A^{-1}(x+(n+1)\alpha,t) \cdots A^{-1}(x-\alpha,t), & n <0
\end{array}\right.    .
$$
Abusing the notation a little bit, for $n\geq 1,$ we define
$$s_n(x,t)=s[A_n(x,t)],~u_n(x,t)=s[A_{-n}(x,t)].$$
We call $s_n$ (respectively, $u_n$) the $n$-step stable (respectively, unstable) direction.

Set $I_0=\R/\Z$ and $g_1(x,t)=s_1(x,t)-u_1(x,t)=\tan^{-1}[t-v(x)]$ for all $t\in J$.  Let $\{\frac{p_n}{q_n}\}_{n\geq1}$ be the continued fraction approximants of $\alpha.$ Fix a  large $N=N(v). $ \cite{wz1} proved the following conclusion by induction. Assume that for $i\geq 1,$  the following objects are well defined:
\begin{enumerate}
\item $i$-th step $critical~points:$
$$C_i(t)=\{c_{i,1}(t),c_{i,2}(t)\}$$
with $c_{i,j}(t)\in I_{i-1,j}(t)$ minimizing $\{|g_i(x,t)|,~x\in I_{i-1,j}(t)\}.$
\item $i$-th step $critical~interval:$
$$I_{i,j}(t)=\{x:|x-c_{i,j}(t)|\leq \frac{1}{2^iq_{N+i-1}^{2\tau}}\}~and~I_i(t)=I_{i,1}(t)\cup I_{i,2}(t).$$
\item $i$-th step $return~times:$
$$q_{N+i-1}\leq r_i^{\pm}(x,t):I_i(t)\rightarrow \Z^+$$
are the first return times (back to $I_i(t)$) after time $q_{N+i-1}-1.$ Here $r_i^+(x,t)$ is the forward return time and $r_i^-(x,t)$ is backward. Let $r_i(t)=\min\left\{r_i^+(t),r_i^-(t)\right\}$ with $r_i^{\pm}(t)=\min_{x\in I_i(t)} r_i^{\pm}(x,t).$
\item ($i+1$)-th step angle $g_{i+1}:$
$$g_{i+1}(x,t)=s_{r_i(t)}(x,t)-u_{r_i(t)}(x,t):D_i\rightarrow \R\mathbb{P}^1,$$
where
$$D_i:=\{(x,t):x\in I_i(t),t\in J\}.$$
\end{enumerate}

The next theorem, which is from \cite{wz1}'s induction theorem, gives the precise description of the several important quantities mentioned above.

\begin{Theorem}[Theorem 3 of \cite{wz1}]\label{theorem12}Given $\alpha\in DC(\kappa,\tau)$, $\e>0$ and a $C^2$ cosine type potential $v$, there exists $\lambda_0(\alpha,v,\e)$ such that  the following holds for $\lambda>\lambda_0$.
\begin{enumerate}
  \item For each $i\geq 1$ and $t\in J,$ it holds that
            \begin{equation}\label{cwc} |c_{i,j}(t)-c_{i+1,j}(t)|<C\lambda^{-\frac{3}{4}r_{i-1}(t)},j=1,2;\end{equation}
\item For each $i\geq 1$, $t\in J,$ and  all $x\in I_{i}(t)$, it holds that
            \begin{equation}\label{mzz}\|A_{\pm r_{i}^{\pm}(x,t)}(x,t)\|>\lambda^{(1-\e)r_{i}^{\pm}(x,t)}\geq\lambda^{(1-\e)q_{N+i-1}}.
            \end{equation}
\item For each $i\geq 1$, $t\in J,$ and  all $x\in I_{i,j}(t)$, it holds that
            \begin{equation}\label{new}
|g_i(x,t)|\geq c|x-c_{i,j}(t)|^{3}\ \ j=1,2.
            \end{equation}
\item  For each $i\geq 1$ and $t\in J$. If $|c_{i,1}(t)-c_{i,2}(t)-k\alpha|\geq \frac{1}{q^{2\tau}_{N+i-1}}$ for all $|k|\leq q_{N+i-1}$, then
\begin{equation}\label{jiaoducha} \|g_{i+1}(\cdot,t)-g_{i}(\cdot,t)\|_{C^2}\leq C\lambda^{-\frac{3}{2}r_{i-1}(t)}.\end{equation}
            \end{enumerate}
\end{Theorem}

Theorem \ref{theorem12} is a simplified version of Theorem 3 in \cite{wz1}. See (57)-(59) in Theorem 3 and Lemma 6 of \cite{wz1}.

\section{Construction of eigenfunctions}
In this section, we construct sufficiently many ``good"  eigenfunctions of $H_{\lambda v,\alpha,\theta}$ by Theorem \ref{theorem12}. We denote by $\Sigma_{\lambda v,\alpha}$ the spectral set of $H_{\lambda v,\alpha,\theta}$ (It does not depend on $\theta$ since $\alpha$ is irrational).
\subsection{The critical points and growth of the transfer matrix}

\begin{Theorem}
Let $\alpha\in DC(\kappa,\tau)$ and $v$ be an even $C^2$ cosine type potential. For any $\e>0$, there exists $\lambda_0(\e,\alpha,v)$ such that if $\lambda>\lambda_0$, and $t\in\lambda^{-1}\Sigma_{\lambda v,\alpha}$, then there exists a strictly increasing continuous surjection
\begin{align}\label{C2}
&\mbox{$c_{\infty}(t):\lambda^{-1}\Sigma_{\lambda v,\alpha}\rightarrow [0,1/2]$},
\end{align}
 and there exist $s_{\infty}(c_{\infty}(t),t), u_{\infty}(c_{\infty}(t),t)\in \R{\mathbb{P}}^1$  if $c_{\infty}(t)\in \Theta_\gamma^\tau$  with
\begin{align}\label{C1}
s_{\infty}(c_{\infty}(t),t)=u_{\infty}(c_{\infty}(t),t),
\end{align}
such that \begin{align}\label{C3}
\|A_n(c_{\infty}(t),t)\hat{s}_{\infty}(c_{\infty}(t),t)\|\geq c\lambda^{-(1-\e)|n|}, \ \ \forall n\in\Z,
\end{align} where $c=c(\kappa,\gamma,\tau,v,
\e)>0$ and $\hat{s}_{\infty}(c_{\infty}(t),t)$ is the unit vector in the direction $s_{\infty}(c_{\infty}(t),t)$.
\end{Theorem}
\begin{pf}
By \eqref{cwc} in Theorem \ref{theorem12} (See also Theorem 3 \cite{wz1} and Theorem 2 in \cite{wz2}), for any $t\in \lambda^{-1}\Sigma_{\lambda v,\alpha}$, for all $j\geq 1$ and $m=1,2$, we have
\begin{align}\label{C4}
|c_{j+1,m}(t)-c_{j,m}(t)|\leq C\lambda^{-\frac{3}{4}r_{j-1}(t)}\leq C\lambda^{-\frac{1}{10}q_{N+j-2}},
\end{align}
thus there exists $c_{\infty,m}(t)$ such that
\begin{align*}
c_{\infty,m}(t)=\lim\limits_{n\rightarrow\infty} c_{n,m}(t).
\end{align*}

 We first prove $c_{\infty,1}(t)=-c_{\infty,2}(t)$, this is because of \eqref{mzz} in Theorem \ref{theorem12}, we have
$$
\|A_{\pm r_j(t)}(c_{j,m}(t),t)\|\geq \lambda^{(1-\e)r_j(t)},
$$
this implies
$$
\|A_{r_j(t)}(c_{j,m}(t),t)\cdot s_{r_j(t)}(c_{j,m}(t),t)\|\leq \lambda^{-(1-\e)r_j(t)},
$$
$$
\|A_{-r_j(t)}(c_{j,m}(t),t)\cdot u_{r_j(t)}(c_{j,m}(t),t)\|\leq \lambda^{-(1-\e)r_j(t)}.
$$
Since $v$ is even, we have
\begin{align*}
&\ \ \ \ \|A_{ r_j(t)}(-c_{j,m}(t),t)\cdot\left(\frac{\pi}{2}-u_{r_j(t)}(c_{j,m}(t),t)\right)\|\\
&=\|A_{-r_j(t)}(c_{j,m}(t),t)\cdot u_{r_j(t)}(c_{j,m}(t),t)\|\leq \lambda^{-(1-\e)r_j(t)},
\end{align*}
this implies that
\begin{equation}\label{use1}
\left|\frac{\pi}{2}-u_{r_j(t)}(c_{j,m}(t),t)-s_{r_j(t)}(-c_{j,m}(t),t)\right|\leq \lambda^{-\frac{3}{2}r_j(t)},
\end{equation}
similarly
\begin{equation}\label{use2}
\left|\frac{\pi}{2}-s_{r_j(t)}(c_{j,m}(t),t)-u_{r_j(t)}(-c_{j,m}(t),t)\right|\leq \lambda^{-\frac{3}{2}r_j(t)}.
\end{equation}
\eqref{use1} and \eqref{use2} imply
$$
|g_{r_j(t)}(c_{j,m}(t),t)-g_{r_j(t)}(-c_{j,m}(t),t)|\leq 2\lambda^{-\frac{3}{2}r_j(t)}.
$$
By \eqref{new} in Theorem \ref{theorem12},
$$
c|c_{j,m}(t)-(-c_{j,m}(t))|^3\leq |g_{r_j(t)}(c_{j,m}(t),t)-g_{r_j(t)}(-c_{j,m}(t),t)|.
$$
Hence
\begin{equation}\label{use3}
|c_{j,m}(t)-(-c_{j,m}(t))|\leq C\lambda^{-\frac{1}{2}r_j(t)}.
\end{equation}
By \eqref{C4}, we have
$$
c_{\infty,1}(t)=-c_{\infty,2}(t).
$$
We  simply denote $c_\infty(t)=c_{\infty,1}(t)$. By the induction theorem (Theorem \ref{theorem12}) we have
$$
c_{\infty}(t)=\begin{cases}
0 &\text{$t=\inf {\lambda^{-1}\Sigma_{\lambda v,\alpha}}$},\\
\frac{1}{2} &\text{$t=\sup {\lambda^{-1}\Sigma_{\lambda v,\alpha}}$}.
\end{cases}
$$
Note that $c_{\infty}(t)$ is continuous on $\lambda^{-1}\Sigma_{\lambda v,\alpha}$, since $c_{j,1}(t)$ converges uniformly to $c_{\infty}(t)$ on $\lambda^{-1}\Sigma_{\lambda v,\alpha}$ and $c_{j,1}$ is continuous, combine these together, we get that $c_{\infty}(t)$ is a continuous surjection from $\lambda^{-1}\Sigma_{\lambda v,\alpha}$ to $[0,\frac{1}{2}]$.

 Now we prove $c_{\infty}(t)$ is increasing on $\lambda^{-1}\Sigma_{\lambda v,\alpha}$, we need the following result: for $t\in \lambda^{-1}\Sigma_{\lambda v,\alpha}$ with $c_{\infty}(t)\in \Theta_\gamma^\tau$, i.e.,
$$
|2c_{\infty}(t)-k\alpha|\geq \frac{\gamma}{(|k|+1)^{\tau}},\ \  \forall k\in\Z,
$$
by \eqref{C4},  there exists $j_0(\gamma)$, such that for all $j>j_0$
\begin{align}\label{C5}
|c_{j,1}(t)-c_{j,2}(t)-k\alpha|\geq \frac{1}{q^{2\tau}_{N+j-1}}, \ \ |k|\leq q_{N+j-1},
\end{align}
thus by Corollary 3 in \cite{wz2},
\begin{equation}\label{use4}
\frac{d \left(c_{j,1}(t)-c_{j,2}(t)\right)}{dt}\geq c>0, \footnote{In case of \eqref{C5}, $\rho_j(t)=c_{j,1}(t)-c_{j,2}(t)$ in \cite{wz2}.}
\end{equation}
in a small neighborhood of $t$. In view of \eqref{use3} and \eqref{use4}, for any sequence $t_1>\cdots>t_n>\cdots$ with $t_n\rightarrow t$, we have $c_{\infty}(t_n)> c_{\infty}(t)$ for all $n$ sufficiently large.

We are now ready to prove the monotonicity of $c_{\infty}(t)$. We prove this by contradiction, otherwise, there exists
$t_1<t_2$ such that $c_{\infty}(t_1)>c_{\infty}(t_2)$, for $\gamma\ll |c_{\infty}(t_1)- c_{\infty}(t_2)|$, there exists $y\in \Theta_\gamma^\tau\cap (c_{\infty}(t_2),c_{\infty}(t_1)) $, let $t'=\sup \{t\in (t_1,t_2): c_{\infty}(t)=y\}$ and for $t\in (t',t_2)$, we have
$$
c_{\infty}(t)<y.
$$
Thus there exists sequence $t_1>\cdots>t_n>\cdots$ with $t_n\rightarrow t'$, such that $c_{\infty}(t_n)< c_{\infty}(t')$ for all $n$ sufficiently large which is a contradiction.

We omit the dependence on $t$ in the following. For any $n\geq q^{100C\tau}_{N+j_0-1}$, let $j_0<j_1<\cdots<j_k=n$ be the return times of $c_{\infty}(t)$ to $I_{N+j_0}$, by \eqref{mzz}, \eqref{new} in Theorem \ref{theorem12} and \eqref{C4}, we have
$$
\|A_n(c_{\infty})\|\geq\prod\limits_{i=0}^{k-1}\|A_{j_i-j_{i-1}}(c_{\infty}+j_{i-1}\alpha)\|q_{N+j_0-1}^{6k\tau}\lambda^{-n+j_{k-1}}\geq c\lambda^{(1-\e)|n|}, \footnote{Here we use the fact that $|n-j_{k-1}|\leq q_{N+j_0-1}^C$, it follows from Lemma 3 of \cite{adz}. }
$$
for some $c=c(\kappa,\gamma,\tau,v,\e)>0$. Thus we have proved \eqref{C3}.  Finally, \eqref{C1} follows from \eqref{jiaoducha} in Theorem \ref{theorem12}.

\end{pf}

\subsection{Construction of  good eigenfunctions}
Recall that $t=\lambda^{-1}E$.
\begin{Definition}\label{good}
For any $\gamma>0$ and $C>0$, a normalized eigenfunction \footnote{We say $u(n)$ is normalized if $\sum_n|u(n)|^2=1$.}  $u(n)$ is said to be $(C,\gamma)$-good, if $$|u(n)|\leq Ce^{-\gamma |n|}$$ for any $n\in\Z$.
\end{Definition}
\begin{Proposition}\label{C2}
Assume $\alpha\in DC(\kappa,\tau)$ and $v$ is an even $C^2$ cosine type potential, for any $\e>0$, there exists $\lambda_0(\e,\alpha,v)$ and $C(\kappa,\gamma,\tau,v,\e)$ such that if $\lambda>\lambda_0$ and $c_{\infty}(t)\in \Theta_\gamma^\tau$, then $H_{\lambda v,\alpha,c_{\infty}(t)}$ has a $(C,(1-\e)\ln\lambda)$-good eigenfunction corresponding to eigenvalue $E=\lambda t$.
\end{Proposition}
\begin{pf}
We denote by $A^{E,\lambda}(\theta)=\begin{pmatrix}E-\lambda v(\theta)&-1\\ 1&0\end{pmatrix}$, for $n\geq 0$, by  \eqref{schcon}, we have
$$
\left|A^{E,\lambda}_{n+1}(c_\infty(t))\cdot \left(B^{-1}(c_\infty(t),t)\cdot s_{n+1}(c_\infty(t),t)\right)\right|\leq C\left\|A_{n+1}(c_\infty(t),t)\right\|^{-1},$$
\begin{align*}
&\left|A^{E,\lambda}_{n+1}(c_\infty(t))\cdot \left(B^{-1}(c_\infty(t),t)\cdot s_{n}(c_\infty(t),t)\right)\right|\\
=&\left|A^{E,\lambda}(c_\infty(t)+n\alpha)A^{E,\lambda}_{n}(c_\infty(t))\cdot \left(B^{-1}(c_\infty(t),t)\cdot s_{n}(c_\infty(t),t)\right)\right|\\
\leq &C\left\|A_n(c_\infty(t),t)\right\|^{-1}.
\end{align*}
Hence
\begin{align*}
&\left|\left(B^{-1}(c_\infty(t),t)\cdot s_{n+1}(c_\infty(t),t)\right)-\left(B^{-1}(c_\infty(t),t)\cdot s_{n}(c_\infty(t),t)\right)\right|\\
\leq& 2C\left\|A_n(c_\infty(t),t)\right\|^{-1}\left\|A_{n+1}(c_\infty(t),t)\right\|^{-1}.
\end{align*}
This implies that
\begin{align*}
&\left|\left(B^{-1}(c_\infty(t),t)\cdot s_{\infty}(c_\infty(t),t)\right)-\left(B^{-1}(c_\infty(t),t)\cdot s_{n}(c_\infty(t),t)\right)\right|\\
\leq &\sum\limits_{k\geq n}2C\left\|A_k(c_\infty(t))\right\|^{-1}\left\|A_{k+1}(c_\infty(t))\right\|^{-1}.
\end{align*}
We have
\begin{align*}
&\left|A^{E,\lambda}_{-n}(c_\infty)\cdot \left(B^{-1}(c_\infty(t),t)\cdot u_{\infty}(c_\infty(t),t)\right)\right|\\
\leq &C\left\|A_{n}(c_\infty(t),t)\right\|^{-1}+\left\|A^{E,\lambda}_{n}(c_\infty(t))\cdot \left(B^{-1}(c_\infty(t),t)\cdot (s_{n}(c_\infty(t),t)-s_{\infty}(c_\infty(t),t))\right)\right\|\\
\leq &C\left(\left\|A_{n}(c_\infty(t),t)\right\|^{-1}+\left\|A_{n}(c_\infty(t),t)\right\|\sum\limits_{k\geq n}\left\|A_k(c_\infty(t),t)\right\|^{-1}\left\|A_{k+1}(c_\infty(t),t)\right\|^{-1}\right).
\end{align*}
By \eqref{C3}, we have
$$
\sum\limits_{k\geq n}\left\|A_k(c_\infty(t),t)\right\|^{-1}\left\|A_{k+1}(c_\infty(t),t)\right\|^{-1}\leq C\lambda^{-2(1-\e)|n|},
$$
for some $C=C(\kappa,\gamma,\tau,v,\e)$.

Similarly,
\begin{align*}
&\left|A^{E,\lambda}_{n}(c_\infty)\cdot \left(B^{-1}(c_\infty(t),t)\cdot s_{\infty}(c_\infty(t),t)\right)\right|\\
\leq &C\left\|A_{-n}(c_\infty(t),t)\right\|^{-1}+\left\|A^{E,\lambda}_{-n}(c_\infty(t))\cdot \left(B^{-1}(c_\infty(t),t)\cdot (u_{n}(c_\infty(t),t)-u_{\infty}(c_\infty(t),t))\right)\right\|\\
\leq &C\left(\left\|A_{-n}(c_\infty(t),t)\right\|^{-1}+\left\|A_{-n}(c_\infty(t),t)\right\|\sum\limits_{k\geq n}\left\|A_{-k}(c_\infty(t),t)\right\|^{-1}\left\|A_{-k-11}(c_\infty(t),t)\right\|^{-1}\right).
\end{align*}
For $E=\lambda t$, we denote by
$$
\begin{pmatrix}
u_E(n+1)\\ u_E(n)
\end{pmatrix}=A_n^{E,\lambda}(c_\infty(t))B^{-1}(c_\infty(t),t)\cdot s_{\infty}(c_\infty(t),t),
$$
then for $n\in\Z$, by \eqref{C1}, we have
\begin{align*}
\left\|\begin{pmatrix}
u_E(n+1)\\ u_E(n)
\end{pmatrix}\right\|&=\|A_n^{E,\lambda}(c_\infty(t))B^{-1}(c_\infty(t),t)\cdot s_{\infty}(c_\infty(t),t)\|\\
&\leq C\lambda^{-(1-\e)|n|}.
\end{align*}
Thus $(u_E)_n$ is a $(C,(1-\e)\ln\lambda)$-good eigenfunction for $H_{\lambda v,\alpha,c_{\infty}(t)}$.
\end{pf}

\section{Completeness arguments}
\subsection{$\mathcal{L}$-measure}In  \cite{gy}, the authors introduced $\mathcal{R}$-measure to prove the arithmetic version of Anderson localization.  Inspired by the idea, we introduce similarly a measure by Proposition \ref{C2} in stead of reducibility in \cite{gy}. The measure, we  call it   $\mathcal{L}$-measure, will play an important role in the proof of arithmetic version of Anderson localization. We next define $E:\T\rightarrow \Sigma_{\lambda v,\alpha}$ as the following:
$$
E(\theta)=\begin{cases}
\lambda c_\infty^{-1}(\theta)&\theta\in [0,\frac{1}{2}],\\
\lambda c_\infty^{-1}(1-\theta)& \theta\in (\frac{1}{2},1].
\end{cases}
$$
Since $c_\infty$ is increasing in the spectrum, $E(\theta)$ takes one value

For every $\tau>1$ and $\gamma>0$, we define  $\mathcal{E}^\tau_\gamma=E(\Theta_\gamma^\tau)$. For any $E\in \mathcal{E}_\gamma^\tau$, we define a vector-valued function $u_E:\mathcal{E}^\tau_\gamma\rightarrow \ell^2(\Z)$ as the following,
\begin{equation}\label{def1}
u_E(n)=\frac{v_E(n)}{\|v_E\|_{L^2}},
\end{equation}
where $v^E$ is the eigenfunction of $H_{\lambda v,\alpha,c_{\infty}(t)}$ constructed in  Proposition \ref{C2}.

For any fixed $\theta\in \Theta^\tau=\cup_{\gamma>0}\Theta^\tau_\gamma$, we denote by $E_m(\theta)=\lambda c_{\infty}^{-1}(T^m\theta)$. We can define the following $\mathcal{L}$-measure,
\begin{Definition}[$\mathcal{L}$-measure]\label{Rm}
\label{defR}
$\nu_{\theta}:\mathcal{B}\rightarrow \R$ is defined as:
$$
\nu_{\theta}(B)=\sum\limits_{m\in N_\theta^B}\frac{|u_{E_m(\theta)}(m)|^2+|u_{E_m(\theta)}(m+1)|^2}{2},
$$
for all $B$ in the Borel $\sigma$-algebra $\mathcal{B}$ of $\R$, where $N_\theta^B=\{m| E_m(\theta)\in B\}$.
\end{Definition}
The $\mathcal{L}$-measure is well defined since the eigenvalues are simple and has the following property.
\begin{Lemma}\label{property} For $a.e.$ $\theta$,
$$\nu_{\theta}(\mathcal{E}^\tau_\gamma)\geq|\Theta_\gamma^\tau|$$ where $|\cdot|$ is the Lebesgue measure.
\end{Lemma}
\begin{pf}
Note that $\nu_{\theta}(\mathcal{E}^\tau_\gamma)$ is measurable in $\theta$ and $\nu_{\theta}(\mathcal{E}^\tau_\gamma)=\nu_{\theta+\alpha}(\mathcal{E}^\tau_\gamma)$. Thus $\nu_{\theta}(\mathcal{E}^\tau_\gamma)=C$ for a.e. $\theta$.

For any $\theta\in \T$, let $N_\theta=\{m| E_m(\theta)\in \mathcal{E}_\gamma^\tau\}$. For any $m\in\Z$, if $m\notin \mathcal{N}_\theta$, let $P_m(\theta)=0$. If $m\in \mathcal{N}_\theta$, let $P_m(\theta)$ be the spectral projection of $H_{\lambda v,\alpha,\theta}$ onto the eigenspace corresponding to $E_{m}(\theta)$. By the definition of $E_m(\theta)$, $u_{E_m(\theta)}(n)$ is an normalized eigenfunction of  $H_{\lambda v,\alpha,T^m\theta}$, thus $T_{-m}u_{E_m(\theta)}(n)$ \footnote{$T_{-m}$ is a translation defined by $T_{-m}u(n):=u(n+m)$.} is an normalized eigenfunction of $H_{\lambda v,\alpha,\theta}$. Now we define a projection operator for any $\theta\in\T$,
$$
P(\theta)=\sum_{T^m\theta\in \Theta_\gamma^\tau}P_m(\theta).
$$
Note that  all these $E$'s in $\mathcal{E}_\gamma^\tau$ are different and  all $P_m(\theta)$ are mutually orthogonal. It follows that  $P(\theta)$ is a projection. Moreover, we have
\begin{align*}
\int_{\T}\frac{\langle P(\theta)\delta_0, \delta_0\rangle+\langle P(\theta)\delta_0, \delta_0\rangle}{2} d\theta&=\int_{\T}\sum\limits_{T^m\theta\in \Theta_\gamma^\tau}\frac{\langle P_m(\theta)\delta_0, \delta_0\rangle+\langle P_m(\theta)\delta_0, \delta_0\rangle}{2} d\theta.
\end{align*}
By Fubini theorem, we have
\begin{align*}
&\int_{\T}\sum\limits_{T^m\theta\in \Theta_\gamma^\tau}\frac{\langle P_m(\theta)\delta_0, \delta_0\rangle+\langle P_m(\theta)\delta_1, \delta_1\rangle}{2} d\theta\\
=&\int_{ \Theta_\gamma^\tau}\sum\limits_{m\in\Z}\frac{\langle P_mT^{-m}(\theta)\delta_0, \delta_0\rangle+\langle P_m(T^{-m}\theta)\delta_1, \delta_1\rangle}{2} d\theta.
\end{align*}
Since $T_{m}H_{\lambda v,\alpha,T^{-m}\theta}T_{-m}=H_{\lambda v,\alpha,\theta}$, we have
\begin{align*}
H_{\lambda v,\alpha,T^{-m}\theta}T_{-m}u_{E_m(\theta)}&=T_{-m}H_{\lambda v,\alpha,\theta}u_{E_m(\theta)}=E_m(\theta)T_{-m}u_{E_m(\theta)}\\
&=E_{m}(T^{-m}\theta)T_{-m}u_{E_m(\theta)}.
\end{align*}
It follows that $T_{-m}u_{E_m}(\theta)$ belongs to the range of $P_{m}(T^{-m}\theta)$, and for each $\delta_n\in\ell^2(\Z)$, we have
$$
\langle P_m(T^{-m}\theta)\delta_n,\delta_n\rangle\geq |\langle T_{-m}u_{E_m(\theta)},\delta_n\rangle|^2.
$$
This implies that
\begin{align*}
&\sum\limits_{m\in\Z}\int_{ \Theta_\gamma^\tau}\frac{\langle P_mT^{-m}(\theta)\delta_0, \delta_0\rangle+\langle P_m(T^{-m}\theta)\delta_1, \delta_1\rangle}{2} d\theta\\
\geq&\sum\limits_{m\in\Z}\int_{\Theta_\gamma^\tau}\frac{|\langle T_{-m}u_{E_m(\theta)}, \delta_0\rangle|^2+|\langle T_{-m}u_{E_m(\theta)}, \delta_1\rangle|^2 }{2} d\theta\\
=&\sum\limits_{m\in\Z}\int_{\Theta_\gamma^\tau}\frac{|\langle T_{m}u_{E_m(\theta)}, \delta_0\rangle|^2+|\langle T_{-m}u_{E_m(\theta)}, \delta_1\rangle|^2 }{2} d\theta.
\end{align*}
Since $u_{E_m(\theta)}$ is a normalized eigenfunction, i.e.,
$$
\sum\limits_{m\in\Z}|\langle T_mu_{E_m(\theta)}, \delta_0\rangle|^2=\sum\limits_{m\in\Z}|\langle T_mu_{E_m(\theta)}, \delta_1\rangle|^2=1.
$$
Hence we have
\begin{align*}
&\int_{\T}\frac{\langle P(\theta)\delta_0, \delta_0\rangle+\langle P(\theta)\delta_1, \delta_1\rangle}{2} d\theta\\
\geq &\sum\limits_{m\in\Z}\int_{\Theta_\gamma^\tau}\frac{|\langle T_{m}u_{E_m(\theta)}, \delta_0\rangle|^2+|\langle T_{-m}u_{E_m(\theta)}, \delta_1\rangle|^2 }{2} d\theta\\
=&|\Theta_\gamma^\tau|.
\end{align*}
Thus
$$
\nu_{\theta}(\mathcal{E}_\gamma^\tau)\geq |\Theta_\gamma^\tau|
$$
for $a.e.$ $\theta$. This finishes  the proof.
\end{pf}
\subsection{Arithmetic version of Anderson localization}
\begin{Lemma}\label{tail}
For any $\epsilon>0$, there exists $N_0(\gamma,\tau,v,\alpha,\epsilon)>0$ such that for all $\theta\in \Theta_\gamma^\tau$, we have
$$
\mathcal{R}_{N_0}\nu_{\theta}(\mathcal{E}^\tau_\gamma):=\sum\limits_{|m|>N_0:T^m\theta\in\Theta_\gamma^\tau}\frac{|u_{E_m(\theta)}(m)|^2+|u_{E_m(\theta)}(m+1)|^2}{2}\leq \epsilon.
$$
\end{Lemma}
\begin{pf}
Note that for any $T^m\theta\in\Theta_\gamma^\tau$, by Proposition \ref{C2},
$$
\frac{|u_{E_m(\theta)}(m)|^2+|u_{E_m(\theta)}(m+1)|^2}{2}\leq Ce^{-\frac{\ln\lambda}{2}|m|}.
$$
Thus for any $\epsilon>0$, there exists $N_0(\gamma,\tau,v,\alpha,\epsilon)>0$ such that for all $\theta\in \Theta_\gamma^\tau$,
$$
\sum\limits_{|m|>N_0:T^m\theta\in\Theta_\gamma^\tau}\frac{|u_{E_m(\theta)}(m)|^2+|u_{E_m(\theta)}(m+1)|^2}{2}\leq \epsilon.
$$
\end{pf}
We define $\mathcal{E}^\tau=\bigcup_{\gamma>0}\mathcal{E}_\gamma^\tau$. Then
\begin{Lemma}\label{continuity}
For any $N>0$ and $\epsilon>0$, there exists $\delta(\gamma,\tau,v,\alpha,N,\epsilon)>0$ such that
$$
\left|\mathcal{T}_N\nu_{\theta}(\mathcal{E}^\tau)-\mathcal{T}_N\nu_{\theta'}(\mathcal{E}^\tau)\right|\leq \epsilon
$$
for any $\theta,\theta'\in\Theta^\tau_\gamma$ with $|\theta-\theta'|\leq \delta$ where $\mathcal{T}_N\nu_{\theta}(\mathcal{E}^\tau):=\nu_{\theta}(\mathcal{E}^\tau)-\mathcal{R}_N\nu_{\theta}(\mathcal{E}^\tau)$.
\end{Lemma}
\begin{pf}
For any fixed $N$ and any $\theta\in\Theta^\tau_\gamma$, we have
$$
\|\theta+k\alpha+n\alpha\|_{\R/\Z}\geq \frac{\gamma}{(|k+n|+1)^\tau}\geq \frac{\gamma(1+N)^{-\tau}}{(|n|+1)^\tau},\ \  |k|\leq N.
$$
Thus $T^k\theta\in \Theta^\tau$ for $|k|\leq N$. By Proposition \ref{theorem12}, we have
$$
\|v_{E_m(\theta)}\|_{\ell^2}\leq C(\gamma,\tau,v,\alpha).
$$
By telescoping and Lemma \ref{tail},  there exists $\delta(\gamma,\tau,v,\alpha,N,\epsilon)$, such that if $|\theta-\theta'|<\delta$, then
\begin{align}\label{redest10}
\left|\|v_{E_k(\theta)}\|^2_{\ell^2}-\|v_{E_k(\theta')}\|^2_{\ell^2}\right|\leq  \frac{\epsilon C^{-4}}{500(2N+1)}.
\end{align}
\begin{align}\label{redest11}
\left|v_{E_k(\theta)}(k)-v_{E_k(\theta')}(k)\right|\leq \frac{\epsilon C^{-4}}{500(2N+1)}, \ \ \forall |k|\leq 2N+1.
\end{align}
\eqref{redest10} and \eqref{redest11} imply for any $|k|\leq N$,
\begin{align*}
|u_{E_k(\theta)}(k)-u_{E_k(\theta')}(k)|&=\left|\frac{v_{E_k(\theta)}(k)}{\|v_{E_k(\theta)}\|_{\ell^2}}-\frac{v_{E_k(\theta')}(k)}{\|v_{E_k(\theta')}\|_{\ell^2}}\right|\\
&=\frac{\left|v_{E_k(\theta)}(k)\|v_{E_k(\theta')}\|_{\ell^2}-v_{E_k(\theta')}(k)\|v_{E_k(\theta)}\|_{\ell^2}\right|}{\|v_{E_k(\theta)}\|_{\ell^2}\|v_{E_k(\theta')}\|_{\ell^2}}\\
&\leq \frac{\epsilon C^{-4}}{500(2N+1)}\frac{\|v_{E_k(\theta)}\|_{\ell^2}+C}{\|v_{E_k(\theta)}\|_{\ell^2}\|v_{E_k(\theta')}\|_{\ell^2}}.
\end{align*}
Thus we have
\begin{align*}
|u_{E_k(\theta)}(k)-u_{E_k(\theta')}(k)|\leq \frac{\epsilon}{100(2N+1)}.
\end{align*}
By Definition \ref{defR} and \eqref{def1}, one has
\begin{align*}
&\ \ \ \ |\mathcal{T}_N\nu_{\theta}(\mathcal{E}^\tau)-\mathcal{T}_N\nu_{\theta'}(\mathcal{E}^\tau)|\\
&=\left|\sum\limits_{|k|\leq N}|u_{E(T^k\theta)}(k)|^2-\sum\limits_{|k|\leq N}|u_{E(T^k\theta')}(k)|^2\right|\\&
+\left|\sum\limits_{|k|\leq N}|u_{E(T^k\theta)}(k+1)|^2-\sum\limits_{|k|\leq N}|u_{E(T^k\theta')}(k+1)|^2\right|\\
&\leq \frac{\epsilon}{50(2N+1)}(2N+1)\leq \epsilon.
\end{align*}
\end{pf}
\begin{Lemma}\label{diophantine}
For any $\theta\in \Theta_\gamma^\tau$, we have $\theta$ is $\Theta^{100\tau}_{\gamma/100}$- homogenous, i.e., $
|(\theta-\sigma,\theta+\sigma)\cap \Theta^{100\tau}_{\gamma/100}|\geq \sigma$ for any $\sigma>0$.
\end{Lemma}
\begin{pf}
Let $\Theta_k=\{\theta\in[0,1):\|2\theta+k\alpha\|_{\R/\Z}<\frac{\gamma}{100(|k|+1)^{100\tau}}\}$, then
$\Theta^{100\tau}_{\gamma/100}=[0,1)\backslash\cup_{k\in\Z}\Theta_k$. Thus for any $\sigma>0$ sufficiently small, we have
$$
(\theta-\sigma,\theta+\sigma)\cap \Theta^{100\tau}_{\gamma/100}=(\theta-\sigma,\theta+\sigma)\backslash\cup_{k\in\Z}\Theta_k.
$$
Notice that if $\Theta_{k}\cap (\theta-\sigma,\theta+\sigma)\neq \emptyset$, then
$$
\frac{\gamma}{(|k|+1)^\tau}\leq\|k\alpha+2\theta\|_{\R/\Z}\leq 4\sigma.
$$
It follows that $|k|\geq (\frac{\gamma}{4\sigma})^{\frac{1}{\tau}}$, thus $|\Theta_k|\leq C(\gamma,\tau)\sigma^{100}$ which implies that
$$
|(\theta-\sigma,\theta+\sigma)\cap \Theta^{100\tau}_{\gamma/100}|\geq 2\sigma-\sum\limits_{k\geq (\frac{\gamma}{\sigma})^{\frac{1}{\tau}}}\Theta_k\geq \sigma.
$$
\end{pf}

\noindent
{\bf Proof of Theorem \ref{main}}: For any fixed $\theta\in \Theta_\gamma^\tau$, by Lemma \ref{diophantine} and Lemma \ref{property}, for any $\gamma_1<\gamma/100$, there exists a sequence $\theta_k\in\Theta^{100\tau}_{\gamma_1}$ such that $\theta_k\rightarrow\theta$ and
$$
\nu_{\theta_k}(\mathcal{E}^{100\tau}_{\gamma_1})\geq\left|\Theta^{100\tau}_{\gamma_1}\right|.
$$
By Lemma \ref{tail}, there exists $N_0(\gamma_1,\tau,,\alpha)>0$ such that
$$
\mathcal{R}_{N_0}\nu_{\theta_k}(\mathcal{E}^{100\tau}_{\gamma_1})\leq \gamma_1.
$$
By Lemma \ref{continuity},
$$
\mathcal{T}_{N_0}\nu_{\theta}(\mathcal{E}^{100\tau})=\lim\limits_{k\rightarrow \infty}\mathcal{T}_{N_0}\nu_{\theta_k}(\mathcal{E}^{100\tau}).
$$
Thus
\begin{align*}
\nu_{\theta}(\mathcal{E}^{100\tau})&\geq \mathcal{T}_{N_0}\nu_{\theta}(\mathcal{E}^{100\tau})\geq \limsup\limits_{k\rightarrow\infty}\mathcal{T}_{N_0}\nu_{\theta_k}(\mathcal{E}^{100\tau}_{\gamma_1})\\
&\geq \left|\Theta^{100\tau}_{\gamma_1}\right|-\gamma_1\geq 1-2\gamma_1.
\end{align*}
Let $\gamma_1\rightarrow 0$, we have
$$
1\leq \nu_{\theta}(\mathcal{E}^{100\tau})\leq\mu^{pp}_{\theta}(\mathcal{E}^{100\tau})\leq\mu_{\theta}(\mathcal{E}^{100\tau})\leq1,  \footnote{We refer readers to Lemma 3.2 in \cite{gy} for the proof of the second inequality.}
$$
where $\mu_{\theta}$ is the spectral measure of $H_{\lambda v,\alpha,\theta}$ defined by
$$
\frac{1}{2}\left(\langle\delta_0,\chi_{B}(H_{\lambda v,\alpha,\theta})\delta_0\rangle+\langle\delta_1,\chi_{B}(H_{\lambda v,\alpha,\theta})\delta_1\rangle\right)=\int_{\R}\chi_{B} d\mu_{\theta},
$$
and  $\mu^{pp}_{\theta}$ is pure point piece of $\mu_{\theta}$. It follows  that $\mu_{\theta}=\mu_{\theta}^{pp}$ for any $\theta\in \Theta_\gamma^\tau$. Thus we finish the proof. \qed

\section*{Acknowledgement}
 J. You was partially supported by National Key R\&D Program of China (2020YFA0713300), NNSF of China (11871286). L. Ge and X. Zhao were partially supported by NSF DMS-1901462. L. Ge was partially supported by AMS-Simons Travel Grant 2020-2022.


\begin{thebibliography}{99}
\bibitem{a} M. Aizenman. Localization at weak disorder: some elementary bounds. {\it Rev. Math. Phys.} {\bf6} (1994), 1163-1182.
\bibitem{ag} M. Aizenman and G.M. Graf. Localization bounds for an electron gas. {\it J. Phys. A: Math. Gen.} {\bf31} (1998), 6783-6806.
\bibitem{am}  M. Aizenman and S. Molchanov. Localization at large disorder and at extreme energies : an elementary derivations. {\it Commun. Math. Phys.} {\bf157} (1993), 245-278.
\bibitem{asfh} M. Aizenman, J. Schenker, R. Friedrich and D. Hundertmark. Finite-volume fractional-moment criteria for Anderson localization. {\it Commun. Math. Phys.} {\bf224} (2001), 219-253.


\bibitem{anderson} P. Anderson. Absence of diffusion in certain random lattices. {\it Phys. Rev.} {\bf109} (1958), 1492-1505.

\bibitem{adz} A. Avila, D. Damanik and Z. Zhang. Singular density of states measures for subshift and quasi-periodic Schr\"odinger operators. {\it Commun. Math. Phys.} \textbf{330} (2014), 469-498.
\bibitem{aj1} A. Avila and S. Jitomirskaya. Almost localization and almost reducibility. {\it J. Eur. Math. Soc} {\bf12} (2010), 93-131.
\bibitem{ayz} A. Avila, J. You and Q. Zhou. Sharp phase transitions for the almost Mathieu operator. {\it Duke Math. J.} {\bf14} (2017), 2697-2718.
\bibitem{aos} J. E. Avron, D. Osadchy and R. Seiler. A topological look at the quantum Hall effect. \textit{Physics today.} {\bf56} (2003), 38-42.




\bibitem{bg} J. Bourgain and M. Goldstein. On nonperturbative localization with quasi-periodic potential. {\it Ann. of Math.} {\bf152} (2000), 835-879.
\bibitem{bj} J. Bourgain and S. Jitomirskaya. Anderson localization for the band model. {\it Geometric Aspects of Functional Analysis.} {\it Lecture Notes in Math.} Springer, Berlin, {\bf1745} (2000), 67-79.
\bibitem{bj02} J. Bourgain and S. Jitomirskaya. Absolutely continuous spectrum for 1D quasiperiodic operators. {\it Invent. Math.} \text{148} (2002), 453-463.
\bibitem{dks} F. Delyon, H. Kunz and B. Souillard. One-dimensional wave equations in disordered media. {\it J. Phys. A: Math. Gen.} {\bf16} (1983), 25-42.
\bibitem{Eli97} L. H. Eliasson. Discrete one-dimensional quasi-periodic Schr\"{o}dinger operators with pure point spectrum. {\it Acta Math.} \textbf{179} (1997), 153-196.

\bibitem{fv} Y. Forman and T. Vandenboom.  Localization and Cantor spectrum for $C^2$ quasi-periodic discrete Schr\"{o}dinger operators. arXiv:2107.05461v1.
\bibitem{fsw} J. Fr\"{o}hlich, T. Spencer and P. Wittwer. Localization for a class of one dimensional quasi-periodic Schr\"{o}dinger operators. {\it Commun. Math. Phys.} \textbf{132} (1990), 5-25.

\bibitem{gk} L. Ge and I. Kachkovskiy. Ballistic transport for one-dimensional quasiperiodic Schr\"odinger operators. arXiv:2009.02896. To appear in {\it Comm. Pure Appl. Math}.
\bibitem{gy} L. Ge and J. You. Arithmetic version of Anderson localization via reducibility. {\it Geom. Funct. Anal.} {\bf30(5)} (2020), 1370-1401.
\bibitem{gyzh1} L. Ge, J. You and X. Zhao. The arithmetic transition conjecture: new proof and generalization. Preprint.
\bibitem{gyz} L. Ge, J. You and Q. Zhou. Exponential dynamical localization: Criterion and applications. arXiv:1901.04258. To appear in {\it Ann. Sci. Ec. Norm. Super. (4)}.
\bibitem{h}  B.I. Halperin. Quantized Hall conductance, current-carrying edge states, and the existence of extended
states in a two-dimensional disordered potential.  {\it Phys. Rev. B.} {\bf25}  (1982), 2185-2190.
\bibitem{j} S. Jitomirskaya. Almost everything about the almost Mathieu operator, II in Proc. of XI Int. Congress of Math. Physics. {\it Int. Press, Somerville, Mass.} (1995), 373-382.
\bibitem{J} S. Jitomirskaya. Metal-Insulator Transition for the almost Mathieu operator. {\it Ann. of Math.} \textbf{150} (1999), 1159-1175.
\bibitem{jk2}S. Jitomirskaya and I. Kachkovskiy. \textit{$L^2$-reducibility and localization for quasiperiodic operators.} {\it Math. Res. Lett.} {\bf 23} (2016), 431-444.
\bibitem{JLiu} S. Jitomirskaya and W. Liu. Universal hierarchical structure of quasi-periodic eigenfuctions. {\it Ann. of Math.} {\bf187(3)} (2018), 721-776.
\bibitem{JLiu1} S. Jitomirskaya and W. Liu. Universal reflective-hierarchical structure of quasiperiodic eigenfunctions and sharp spectral transition in phase. arXiv:1802.00781.
\bibitem{js} S. Jitormiskya and  B. Simon. Operators with singular continuous spectrum. III. Almost periodic Schr\"odinger operators. {\it Commun. Math. Phys.} {\bf 165} (1994), 201-105.
\bibitem{fk} A. Klein and A. Figotin. Midgap defect modes in dielectric and acoustic media. {\it SIAM. J. Appl. Math.}
{\bf58} (1998), 1748-1773.
\bibitem{K} S. Klein. Anderson localization for the discrete one-dimensional quasi-periodic Schr\"odinger operator with potential defined by a Gevrey-class function.  {\it J. Funct. Anal.} {\bf218} (2005), 255-292.

\bibitem{ks} H. Kunz and B. Souillard. Sur le spectre des op$\acute{e}$rateurs aux diff$\acute{e}$rences finies al$\acute{e}$atoires. {\it Commun. Math. Phys.} {\bf78} (1980/81), 201-246.
\bibitem{ntw} Q. Niu, D.J. Thouless and. Y.S. Wu. Quantized  Hall conductance as a topological invariant. {\it Phys.Rev.
B.} {\bf31}  (1985), 3372-3377.






\bibitem{oa} D. Osadchy and J.E.  Avron. Hofstadter butterfly as quantum phase
diagram. \textit{J. Math Phys.} {\bf42} (2001), 5665-5671.
\bibitem{Sinai} Ya. G. Sinai. Anderson localization for one-dimensional difference Schr\"{o}dinger operator with quasi-periodic potential. {\it J. Stat. Phys.} \textbf{46} (1987), 861-909.
\bibitem{wz1} Y. Wang and Z. Zhang. Uniform positivity and continuity of Lyapunov exponents for a class of $C^2$ quasiperiodic Schr\"odinger cocycles. {\it J. Funct. Anal.} {\bf268} (2015), 2525-2585.
\bibitem{wz2} Y. Wang and Z. Zhang. Cantor spectrum for a class of $C^2$ quasiperiodic Schr\"odinger operators. {\it Int. Math. Res. Not.} {\bf 8} (2017), 2300-2336.
\bibitem{young} L.-S. Young. Lyapunov exponent fir some quasi-periodic cocycles. {\it Ergodic Theory Dynam. Systems} {\bf17} (1997), 483-504.
\end{thebibliography}
\end{document}